\documentclass[smallextended,envcountsect]{svjour3}

\usepackage{amssymb}

\smartqed

\newcommand{\N}{{\mathbb N}}
\newcommand{\R}{{\mathbb R}}
\newcommand{\X}{{\mathbb X}}
\newcommand{\Y}{{\mathbb Y}}
\newcommand{\E}{{\mathbb E}}

\newcommand{\fix}{{\rm fix\,}}
\newcommand{\cl}{{\rm cl\,}}

\newcommand{\gph}{{\rm gph\,}}

\begin{document}

\title{Approximate solutions of quasiequilibrium problems in Banach spaces}

\author{M. Castellani \and M. Giuli}

\institute {M. Castellani \and M. Giuli \at Department of Information Engineering, Computer Science and Mathematics, L'Aquila, Italy \\
\email{massimiliano.giuli@univaq.it}}

\maketitle

\begin{abstract}
In this note we show that a recent existence result on quasiequilibrium problems, which seems to improve deeply some well-known results, is not correct.
We exhibit a counterexample and we furnish a generalization of a lemma about continuous $\varepsilon$-minimizers of quasiconvex functions depending on a parameter. This allows to establish an existence result of approximate solutions of quasiequilibrium problems.
\end{abstract}
\keywords{Quasiequilibrium problems \and Lower semicontinuous set-valued maps \and Fixed points \and Approximate solutions}

\section{Introduction}\label{sect:introduction}

The well-known Ky Fan inequality \cite{Fa72} is a very general mathematical model which embraces the formats of several disciplines.
This general problem was named equilibrium problem by Blum and Oettli \cite{BlOe94}, who stressed this unifying feature and provided a thorough investigation of its theoretical properties.

A quasiequilibrium problem is an equilibrium problem in which the constraint set is subject to modifications depending on the considered point.
This problem setting encompasses many relevant problems as special cases, among which variational and quasivariational inequalities, generalized Nash equilibrium problems, mixed quasivariational-like inequalities and so on.

There exists quite a bit of literature on the equilibrium problems (see \cite{BiCaPaPa13} for a recent survey on existence results and solution methods), but not so many studies on quasiequilibrium problems even if a first existence result was stated, to the best of our knowledge, at 1976 in a seminal paper of Mosco \cite{Mo76}.
This is due to the fact that the set of solutions of a quasiequilibrium problem coincides with the set of fixed points a map
which is tricky to handle.

An interesting existence result for quasiequilibrium problems on $\R^n$ was established in \cite{Cu97} where the upper semicontinuity of the set-value map $K$ is not assumed.
The proof of the theorem depends heavily on the finite dimensionality of the space.
In a recent paper \cite{GePa11}, the authors established an existence result of approximate solutions of quasiequilibrium problem defined on a Banach space and, applying the Maximum Theorem \cite{Be59}, they extended the result in \cite{Cu97} to infinite dimensional spaces.
Unfortunately the first result is incorrect and the mentioned improvement fails.

In the present note a counterexample which contradicts the existence of approximate solutions is exhibited.
The mistake in the proof stems from an incorrect application of a result \cite{Ge05} about continuous $\varepsilon$-minimizers of quasiconvex functions depending on a parameter which seems worthy in itself.
We furnish a slightly improved version of this result which allows to establish a new existence result of approximate solutions of quasiequilibrium problems.

\section{The results}\label{sect:results}

In order to better present the rest of this note, we recall notations and definitions which will be used in the sequel.
A set-valued map $F:\X\rightrightarrows\Y$, where $\X$ and $\Y$ are two topological spaces, is lower semicontinuous if for any open subset $A$ of $\Y$ the lower inverse set
\begin{displaymath}
F^{-1}(A)=\{x\in\X :F(x)\cap A\ne\emptyset\}
\end{displaymath}
is open.
We also introduce the notation $F^{-1}(y)=F^{-1}(\{y\})$ for the lower section of $F$ at $y$ and
$\gph(F)=\{(x,y)\in \X\times \Y: y\in F(x)\}$ for the graph of $F$.
Clearly set-valued maps with open graph have open lower sections which in turn are lower semicontinuous.
The set of fixed points of $F$ will be denoted by $\fix F$.
We denote by $B$ the open unit ball in a metric space and by $\cl$ the closure of a set.
Finally, a function $\varphi:C\rightarrow\R$ is quasiconcave on the convex subset $C$ of $\X$ if all the upper level sets are convex; $\varphi$ is quasiconvex if $-\varphi$ is quasiconcave.

Hereafter, we consider the following quasiequilibrium problem:
\begin{equation}\label{eq:qep}
\mbox{find } x\in K(x) \mbox{ such that } f(x,y) \leq 0 \mbox{ for all } y\in K(x),
\end{equation}
where $K:C\rightrightarrows C$ is a set-valued map defined on a nonempty subset $C$ of a suitable topological space $\E$ and $f:C\times C\rightarrow\R$ is a bifunction.
An interesting existence result for (\ref{eq:qep}) was established in \cite{Cu97} where the author replaced the usual upper semicontinuity and the closed-valuedness assumption on $K$ with the closedness of the set $\fix K$.
\begin{theorem}\label{th:Cu97}
Take $\E=\R^n$ and suppose that $C$ is convex and compact, the set-valued map $K$ is lower semicontinuous with convex and nonempty images and $\fix K$ is closed.
Assume that $f$ has the properties:
\begin{description}
\item[{\rm (i)}] the sublevel set $\{(x,y)\in C\times C:f(x,y)\le 0\}$ is closed;
\item[{\rm (ii)}] the function $f(x,\cdot)$ is quasiconcave on $K(x)$ for every $x\in C$;
\item[{\rm (iii)}] $f(x,x)\le 0$ for all $x\in C$.
\end{description}
Then (\ref{eq:qep}) has solution.
\end{theorem}
The proof of Theorem \ref{th:Cu97} depends heavily on the finite dimensionality of the space.
In a very recent paper \cite{GePa11}, the following result, which seems to improve deeply Theorem \ref{th:Cu97}, was stated.
\begin{theorem}\label{th:GePa11}
Let $\E$ be a Banach space and suppose that $C$ is convex and compact, and the set-valued map $K$ is lower semicontinuous with convex and nonempty images.
Assume that $f$ has the properties:
\begin{description}
\item[{\rm (i)}] the function $f(\cdot,y)$ is lower semicontinuous for every $y\in C$;
\item[{\rm (ii)}] the function $f(x,\cdot)$ is quasiconcave for every $x\in C$;
\item[{\rm (iii)}] there exists $\varepsilon_0>0$ such that $f(x,x)\le 0$ for all $x\in (K(x)+\varepsilon_0B)\cap C$.
\end{description}
Then
\begin{description}
\item[{\rm (a)}] for each $\varepsilon>0$ there exists $x_\varepsilon\in\fix K$ such that
\begin{displaymath}
f(x_\varepsilon,y)\le\varepsilon,\qquad\forall y\in K(x_\varepsilon);
\end{displaymath}
\item[{\rm (b)}] if $\fix K$ is closed and $f$ is lower semicontinuous then (\ref{eq:qep}) has solution.
\end{description}
\end{theorem}
Notice that Assertion (b), whose proof relies on an application of the Maximum Theorem \cite{Be59} to Assertion (a), would be a generalization of Theorem \ref{th:Cu97} to infinite dimensional spaces.
Unfortunately Assertion (a) is incorrect as the following example shows at once and the mentioned improvement fails.
\begin{example}\label{ex}
Take $\E=\R$, $C=[0,4]$, $f(x,y)=y-x$ and let the set-valued map $K$ be defined by
\begin{displaymath}
K(x)=\left\{\begin{array}{ll}
[0,4] & \mbox{if } x\in[0,2)\\
\vphantom{.}[0,1) & \mbox{if } x\in[2,3]\\
\vphantom{.}[0,3) & \mbox{if } x\in (3,4]
\end{array}\right.
\end{displaymath}
Since $f$ is linear it is lower semicontinuous with respect to the first variable and quasiconcave with respect to the second variable.
Besides $f(x,x)=0$ for all $x\in C$.
On the other hand $\gph K$ is open in $C\times C$ and therefore $K$ has open lower sections and it is lower semicontinuous with $\fix K=[0,2)$.
However, for every $x\in\fix K$, the point $y=4$ belongs to $K(x)$ and $f(x,4)>2$.
Hence there are not $\varepsilon$-solutions of the quasiequilibrium problem when $\varepsilon\in (0,2]$.
\end{example}
As claimed before, the oversight in the proof of Theorem \ref{th:GePa11} relies on an incorrect application of \cite[Lemma 2.1]{Ge05}.
In order to get the existence of $\varepsilon$-solutions of the quasiequilibrium problem (\ref{eq:qep}), we prove a more generalized version of this lemma.
\begin{theorem}\label{le:newGe05}
Suppose that $\X$ is a paracompact topological space, $\Y$ is a locally convex Hausdorff topological vector space, $C_\Y$ is a convex subset of $\Y$, $C_\X$ is a closed subset of $\X$, $F:\X\rightrightarrows C_\Y$ is a lower semicontinuous set-valued map with convex nonempty images, $\varepsilon$ and $\varepsilon'$ two lower semicontinuous functions from $\X$ to $(0,+\infty)$,  $V\subseteq\Y$ is a balanced convex open neighbourhood of the origin, and the functions $f:\X\times C_\Y\rightarrow\R$ and $g:\X\rightarrow\R$ satisfy the following conditions:
\begin{description}
\item[{\rm (i)}]  the function $f(x,\cdot)$ is quasiconvex for every $x\in C_\X$;
\item[{\rm (ii)}] the function $f(\cdot,y)$ is upper semicontinuous for every $y\in C_\Y$;
\item[{\rm (iii)}] $g$ is lower semicontinuous and
\begin{equation}\label{eq:pardalos}
g(x)\ge\inf_{y\in F_{\varepsilon,V}(x)\cap C_\Y}f(x,y),\qquad\forall x\in C_\X
\end{equation}
where $F_{\varepsilon,V}(x)=F(x)+\varepsilon(x)V$.
\end{description}
Then
\begin{description}
\item[{\rm (a)}]
there exists a continuous selection $\varphi_{\varepsilon,V}:\X\rightarrow C_\Y$ of $F_{\varepsilon,V}$ such that
\begin{equation}\label{eq:pandorata}
f(x,\varphi_{\varepsilon,V}(x))<g(x)+\varepsilon'(x),\qquad\forall x\in C_\X;
\end{equation}
\item[{\rm (b)}]
if $F$ has open lower sections, then there exists a continuous selection $\varphi_{\varepsilon,V}$ of $F$ satisfying (\ref{eq:pandorata}).
\end{description}
\end{theorem}
{\it Proof }
(a) This assertion is equivalent to show that the set-valued map $H_{\varepsilon,\varepsilon',V}:\X\rightrightarrows\Y$ defined by
\begin{displaymath}
H_{\varepsilon,\varepsilon',V}(x)=
\left\{
\begin{array}{ll}
G_{\varepsilon,\varepsilon',V}(x) & \mbox{ if } x\in C_\X \\
F_{\varepsilon,V}(x)\cap C_\Y& \mbox{ if } x\notin C_\X
\end{array}
\right.
\end{displaymath}
where $G_{\varepsilon,\varepsilon',V}(x)=F_{\varepsilon,V}(x)\cap \{y\in C_\Y:f(x,y)<g(x)+\varepsilon'(x)\}$, has a continuous selection and it will be enough to verify that all the assumptions of the Browder Selection Theorem (see for instance \cite[Theorem 17.63]{AlBo06}) are satisfied.
The nonemptiness of $H_{\varepsilon,\varepsilon',V}(x)$ descends from (\ref{eq:pardalos}) instead $H_{\varepsilon,\varepsilon',V}(x)$ is convex since both $F_{\varepsilon,V}(x)$
and $\{y\in C_\X:f(x,y)<g(x)+\varepsilon'(x)\}$ are convex sets from assumption (i).

It remains to show that $H_{\varepsilon,\varepsilon',V}$ has open lower sections.
Fix $y\in C_\Y$, first we prove that $F^{-1}_{\varepsilon,V}(y)$ is open.
Take $x\in F^{-1}_{\varepsilon,V}(y)$ and choose $\delta,r,s >0$ such that $y\in F(x)+(\varepsilon(x)-\delta)V$ and $r+s=\varepsilon(x)-\delta$.
From the lower semicontinuity of $\varepsilon$ there exists an open neighbourhood $U_x$ of $x$ such that $\varepsilon(x')>\varepsilon(x)-\delta$ for every
$x'\in U_x$.
Since the lower semicontinuity of $F$ implies the lower semicontinuity of the set-valued map $F_{r,V}$ defined by $F_{r,V}(x)=F(x)+rV$,
this guarantees that the set $U_x\cap F^{-1}_{r,V}(y+sV)$ is an open neighbourhood of $x$ which is contained in $F^{-1}_{\varepsilon,V}(y)$.
The lower semicontinuity of $\varepsilon'$ and assumptions (ii) and (iii) guarantee that the set $\{x\in C:f(x,y)<g(x)+\varepsilon'(x)\}$ is open and hence
$G^{-1}_{\varepsilon,V}(y)$ is also open since intersection of two open sets.
Finally take $x\in H^{-1}_{\varepsilon,\varepsilon',V}(y)$ and distinguish two different cases.
If $x\in C_\X$, the open set $G^{-1}_{\varepsilon,\varepsilon',V}(y)$ is contained in $H^{-1}_{\varepsilon,\varepsilon',V}(y)$ which is a neighborhood of $x$.
Otherwise, if $x\notin C_\X$ then $F^{-1}_{\varepsilon,V}(y)\setminus C_\X\subseteq H^{-1}_{\varepsilon,\varepsilon',V}(y)$ and it is an open
neighborhood of $x$.
In both cases, the set $H^{-1}_{\varepsilon,\varepsilon',V}(y)$ is a neighborhood of $x$ and the claim follows.

The assertion (b) is established in a similar way replacing $F_{\varepsilon,V}$ with $F$ in the definition of $H_{\varepsilon,\varepsilon',V}$.
\qed

\begin{remark}
Theorem \ref{le:newGe05} recovers Lemma 2.1 in \cite{Ge05} where the author assumes $\Y$ a Banach space, $C_\X=\X$, $C_\Y$ closed and $\varepsilon=\varepsilon'$.
In regards to Assertion (b), the assumption of openness of the images of $F$ assumed in \cite{Ge05} does not guarantee the openness of the lower sections of $F$ and also the existence of a continuous selection as showed in \cite[Example 6.3]{Mi56}.
Nevertheless this happens when $\Y=\R^n$, since lower semicontinuous maps with open and convex values must have open lower sections.
\end{remark}

Returning back to Theorem \ref{th:GePa11}, the mistake in its proof resides in the fact that the authors apply Lemma \ref{le:newGe05} choosing
$\X=C_\X=C_\Y=C$, $\Y=\E$, $\varepsilon=\varepsilon'$ and $g$ constant with
\begin{displaymath}
g(z)=\sup_{x\in\fix K}\inf_{y\in (K(x)+\varepsilon B)\cap C}(-f(x,y)),\qquad\forall z\in C
\end{displaymath}
but $g$ does not satisfy (\ref{eq:pardalos}).
For instance if we consider $E$, $C$, $f$ and $K$ as in Example \ref{ex}, we have $g(z)=-1$,
\begin{displaymath}
\inf_{y\in (K(x)+(-\varepsilon,\varepsilon))\cap [0,4]}(x-y)=
\left\{\begin{array}{ll}
x-4 & \mbox{if } x\in[0,2)\\
x-1-\varepsilon & \mbox{if } x\in[2,3]\\
x-3-\varepsilon & \mbox{if } x\in(3,4]
\end{array}\right.
\end{displaymath}
and (\ref{eq:pardalos}) fails in  $[2,4]$.

The last part of this note is devoted to establish an existence result by exploiting Theorem \ref{le:newGe05}.
Before giving the result we need the following lemma.
\begin{lemma}\label{le:berge}
Let $\X$ and $\Y$ be two topological spaces.
Assume that $F:\X\rightrightarrows\Y$ has open lower sections
and $f:\X\times\Y\rightarrow \R$ is lower semicontinuous with respect to the first variable.
Then the extended-valued function $m$ defined by
\begin{displaymath}
m(x) = \sup_{y\in F(x)}f(x,y)
\end{displaymath}
is lower semicontinuous.
\end{lemma}
{\it Proof }
Fix $\alpha\in\R$ and $x\in \X$ such that $m(x)>\alpha$, therefore there exists $y\in F(x)$ such that $f(x,y)>\alpha$.
The lower semicontinuity of $f(\cdot,y)$ implies the existence of a neighborhood $U'_x$ of $x$ such that
$f(x',y)>\alpha$ for all $x'\in U'_x$.
On the other hand, since $F$ has open lower sections there exists a neighborhood $U''_x$ of $x$ such that $y\in F(x')$ for all
$x'\in U''_x$.
Choosing $U_x=U'_x\cap U''_x$ we get
\begin{displaymath}
m(x')=\sup_{y'\in F(x')}f(x',y')\ge f(x',y)>\alpha \qquad \forall x'\in U_x
\end{displaymath}
and the proof is complete.
\qed
Now we are in the position to prove our existence result.
\begin{theorem}\label{th:GePa11new}
Let $\E$ be a Banach space and suppose that $C$ is convex and compact, and the set-valued map $K$ is lower semicontinuous with convex and nonempty images.
Assume that $f$ has the properties:
\begin{description}
\item[{\rm (i)}] the function $f(\cdot,y)$ is lower semicontinuous for every $y\in C$;
\item[{\rm (ii)}] the function $f(x,\cdot)$ is quasiconcave for every $x\in C$;
\item[{\rm (iii)}] there exists $\varepsilon_0>0$ such that $f(x,x)\le 0$ for all $x\in (K(x)+\varepsilon_0B)\cap C$.
\end{description}
Then for each $0<\varepsilon<\varepsilon_0$ there exists $x_\varepsilon\in\cl \fix K_{\varepsilon}$ such that
\begin{displaymath}
f(x_\varepsilon,y)\le 0,\qquad\forall y\in K_{\varepsilon}(x_\varepsilon)
\end{displaymath}
where $K_\varepsilon:C\rightrightarrows C$ is defined by $K_{\varepsilon}(x)=(K(x)+\varepsilon B)\cap C$.
\end{theorem}
{\it Proof }
Notice that since $K_{\varepsilon}$ has open lower sections for every $\varepsilon >0$, the Browder fixed point theorem guarantees that $\fix K_{\varepsilon}\neq \emptyset$.
For each fixed $0<\varepsilon <\varepsilon_0$ and $n\in\N$,
apply Assertion (a) of Theorem \ref{le:newGe05} to the function $-f$ with $\X=C=C_\Y$, $\Y=\E$, $C_\X=\cl \fix K_{\varepsilon}$, $F=K$, $V=B$ and
\begin{displaymath}
g(z)=\sup_{x\in\cl \fix K_{\varepsilon}}\inf_{y\in K_\varepsilon(x)}(-f(x,y))=-\inf_{x\in\cl \fix K_{\varepsilon}}\sup_{y\in K_\varepsilon(x)}f(x,y) \qquad \forall z\in C.
\end{displaymath}
From Lemma \ref{le:berge} the function $x\mapsto\sup_{y\in K_\varepsilon(x)}f(x,y)$ is proper and lower semicontinuous, so the infimum is attained and $g(z)$ is finite.
Then there exists a continuous selection $\varphi_{\varepsilon}:C\rightarrow C$ of $K_{\varepsilon}$ such that
\begin{equation}\label{eq:appl}
-f(x,\varphi_{\varepsilon}(x))<g(x)+\frac{1}{n},\qquad\forall x\in \cl \fix K_{\varepsilon}.
\end{equation}
The Schauder fixed point theorem ensures the existence of a point $x=\varphi_{\varepsilon}(x)\in K_{\varepsilon}(x)$.
Evaluating (\ref{eq:appl}) at such fixed point and taking limit as $n\to \infty$, we have
\begin{displaymath}
-\inf_{x\in\cl \fix K_{\varepsilon}}\sup_{y\in K_\varepsilon(x)}f(x,y)\ge -f(x,x)\ge0.
\end{displaymath}
Arguing as before the infimum is attained and the theorem is proved.
\qed

We end showing that, in the setting of Theorem \ref{th:GePa11new}, approximate solutions exist even when there are no solutions.
Let $\E$, $C$, $f$ and $K$ be given as in Example \ref{ex}.
Again $\gph K$ is open in $C\times C$, $f$ is linear, $\fix K=[0,2)$ and (\ref{eq:qep}) has no solutions.
However, for every $0<\varepsilon<1$, $\fix K_\varepsilon=[0,2)\cup(3,3+\varepsilon)$ and the point
$x_\varepsilon=3+\varepsilon$ belongs to $\cl\fix K_\varepsilon$ and verifies $f(x_\varepsilon,y)\le0$ for all
$y\in K_\varepsilon(x_\varepsilon)=[0,3+\varepsilon)$.



\end{document}